\documentclass{amsart}

\usepackage{amsmath, amsthm, amsopn, amsfonts}

\usepackage{graphicx}

\usepackage{graphpap}

\def\pasdegrille{\let\grille = \pasgrille}

\def\aat#1#2#3{

\divide \dimen1 by 48

\dimen3=\dimen1

\multiply \dimen1 by #1

\advance \dimen1 by -\dimen3

\divide \dimen1 by 101

\multiply \dimen1 by 100

\divide \dimen2 by \count11

\multiply \dimen2 by #2

\setbox0=\hbox{#3}\ht0=0pt\dp0=0pt

  \rlap{\kern\dimen1 \vbox to0pt{\kern-\dimen2\box0\vss}}\dimen1= \wd1

\dimen2=\ht1}

\def\pasgrille{

\count12= \dimen1

\divide \count12 by 50

\divide \dimen2 by \count12

\count11 =\dimen2

\

\divide \dimen1 by 48

\setlength{\unitlength}{\dimen1}

\smash{\rlap{\ }}

\dimen1= \wd1

\dimen2=\ht1

}

\def\grille{

\count12= \dimen1

\divide \count12 by 50

\divide \dimen2 by \count12

\count11 =\dimen2

\

\divide \dimen1 by 48

\setlength{\unitlength}{\dimen1}

\smash{\rlap{\graphpaper[1](0,0)(50, \count11)}}

\dimen1= \wd1

\dimen2=\ht1

}

\pasdegrille

\usepackage{amssymb}

\ifx\pdfoutput\undefined

\usepackage{graphicx}  \DeclareGraphicsExtensions{.ps}

\else

\usepackage{amsfonts}

\usepackage[pdftex]{graphicx}  \DeclareGraphicsExtensions{.pdf,.mps}

\fi

\usepackage{epic, eepic}

\def\squarebox#1{\hbox to #1{\hfill\vbox to #1{\vfill}}}

\newcommand{\1}{{\bold 1}}

\newcommand{\CI}{{\mathcal C}^\infty }

\newcommand{\CIc}{{\mathcal C}^\infty_{\rm{c}} }

\newcommand{\SP}{{\mathbb S}}

\newcommand{\rest}{\!\!\restriction}

\theoremstyle{plain}

\newtheorem{thm}{Theorem}

\newtheorem{prop}{Proposition}[section]

\theoremstyle{definition}

\newtheorem{rem}{Remark}

\numberwithin{equation}{section}

\title[Eigenfunctions for partially rectangular billiards]
{Eigenfunctions for partially rectangular billiards}
\author[N. Burq]{Nicolas Burq}
\address{Universit{\'e} Paris Sud,
Math{\'e}matiques,
B{\^a}t 425, 91405 Orsay Cedex}
\email{Nicolas.burq@math.u-psud.fr}
\author[M. Zworski]{Maciej Zworski}
\address{Mathematics Department, University of California \\
Evans Hall, Berkeley, CA 94720, USA}
\email{zworski@math.berkeley.edu}

\usepackage{amssymb}
\usepackage{amsmath, amsthm, amsopn, amsfonts}

\def\11{{\rm 1~\hspace{-1.4ex}l} }

\begin{document}

\maketitle   

\section{Introduction}   
\label{in}
In this note we further develop the idea of using a ``black box''
point of view \cite{BZ2} to study eigenfunctions for billiards
which have rectangular components: they include the Bunimovich 
billiard, the Sinai billiard, and the recently popular pseudointegrable billiards
\cite{Bog}.

\begin{figure}[ht]
\includegraphics[width=10cm]{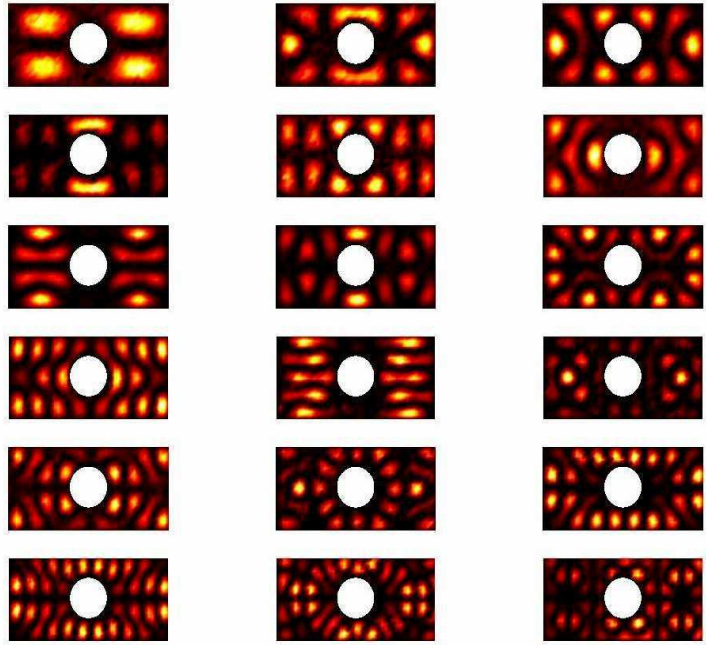}
\caption{Experimental images of eigenfunctions in a Sinai billiard
microwave cavity -- see {\tt http://sagar.physics.neu.edu}. We see
that there is always a non-vanishing presence near the boundary of the
obstacle as predicted by Theorem \ref{t:s} below.}
\label{fig:sin}
\end{figure}

By a partially rectangular billiard we mean a connected planar domain,
$ \Omega $, with a piecewise smooth boundary, which contains a  rectangle,
$ R \subset \Omega $, such that if we decompose the boundary of 
$ R $, into pairs of parallel segments, $ \partial R = \Gamma_1 \cup \Gamma_2 $,
then $ \Gamma_i \subset \partial \Omega $, for at least one $ i$.
Motivated by the general theory of \cite{BZ2} we have used 
elementary methods \cite{BZ3} to show that for such domains 
the eigenfunctions of the Dirichlet, Neumann,
or periodic Laplacian, cannot concentrate in the rectangle, away from the remaining
two sides of the rectangle -- see Theorem \ref{t:1} below. 

In this note we show how a combination of this elementary result with 
the now standard, but highly non-elementary, 
propagation results of Melrose-Sj\"ostrand \cite{MeSj} and Bardos-Lebeau-Rauch \cite{BLR},
gives improved results in some interesting situations. That was already indicated, in
a special case, in \cite[Theorem 3$'$]{BZ2} but here we give an independent
and more general presentation. For the motivation coming from {\em quantum chaos}
we suggest \cite{Do},\cite{Ze},\cite{BZ3}, and references given there.

\medskip
\noindent
{\sc Acknowledgments.} This work was supported in part by 
the National Science Foundation  under the grant  DMS-0200732. We would also
like to thank Srinivas Sridhar for letting us use the experimental 
images shown in Fig.\ref{fig:sin}.  The 
second author is also grateful to Universit\'e de Paris-Sud, Orsay,
for its generous hospitality in September 2003.

\section{Preliminaries}
\label{pr} 
In this section we will recall the basic control result \cite{Bu92},\cite{BZ2}
for rectagles, and the propagation results \cite{MeSj},\cite{BLR},\cite{BuIMRN},\cite{BuGe} for 
billiards. Since in the specific application presented in Sect.\ref{app} we only use propagation
away from the boundary only that, easier, case will be reviewed.

The following result~\cite{Bu92} is related to some earlier control results 
of Haraux~\cite{Ha} and  Jaffard~\cite{Ja}\footnote{We remark that as noted in \cite{Bu92} 
the result holds for 
any product manifold $M= M_{x}\times M_{y}$, and the proof is essentially the same.}
\begin{prop}
\label{p:1} Let $\Delta$ be  the Dirichlet, Neumann, or periodic
Laplace operator on the rectangle $R= [0, 1]_{x} \times [0,a]_{y}$. 
Then for any open non-empty  $\omega\subset R$ of the form $
\omega= \omega_{x} \times [0,a]_{y}$ , there 
exists $C$ such that for any solutions of
\begin{equation}
(\Delta -z) u =f \ \text{ on $R$}, \ u \rest_{\partial R}=0
\end{equation}
we have
\begin{equation}
\label{eq:6.12}\|u\|^2_{{L^2(R)}}\leq C \left(\|f \|^2_{H^{-1} (
[0,1]_{x}; L^2([0,a]_{y})) } +
\|u\rest _{\omega}\|^2_{{L^2(\omega)}} \right)
\end{equation}
\end{prop}
\begin{proof}
We will consider the Dirichlet case (the proof is the same in the other two 
cases) and 
decompose $u,f$ in terms of 
the basis of $L^2([0,a])$ formed by the Dirichlet eigenfunctions
$e_{k}(y)=  { \sqrt {{2}/a}}\sin(2k\pi y/a)$,
\begin{equation}
u(x,y)= \sum_{k}e_{k}(y) u_{k}(x), \qquad f(x,y)= \sum_{k}e_{k}(y) f_{k}(x)
\end{equation}
we get for $u_{k}, f_{k}$ the equation
\begin{equation}\label{estres.1}
\left(\Delta_{x}-\left(z+\left({2k\pi}/{a}\right)^2\right)\right)u_{k}= f_{k},\qquad u_{k}(0)=u_k(1)=0
\end{equation}
We now claim that 
\begin{equation}
\label{eq:cont}
\|u_{k}\|^2_{{L^2([0,1]_{x})}}\leq C \left(\|f_k \|^2_{H^{-1}([0,1]_{x})} +
\|u_k \rest _{\omega_{x}}\|^2_{{L^2(\omega)}}\right) 
\end{equation}
from which, by summing the squares in $k$, we get~\eqref{eq:6.12}.

To see \eqref{eq:cont} we can use the propagation result below in dimension
one, but in this case an elementary calculation is easily available -- see
\cite{BZ3}.
\end{proof}

To state the propagation theorem in the form sufficient for our applications
we follow \cite{BuIMRN} and introduce microlocal defect measures.

Consider for  $a(x, \xi)\in \CIc ( {\mathbb R}^{ 2d})$ and $\varphi\in 
 \CIc ( {\mathbb R}^d)$ equal to $1$ near the $x$-projection of the support of $a$.
To the symbol $a $ we associate 
the family of operators 
 $\text{ Op}_{ \varphi}(a)(x, hD_{x})$ defined by
\begin{equation}
\label{eq2.6}\text{Op}_{\varphi}(a)(x, hD_{x})f= \frac{ 1}{ (2 \pi)^d }\int e^{ix\cdot \xi}a(x, h\xi)\widehat {\varphi f}(\xi) d\xi
\end{equation}
By the symbolic calculus the operator $\text{ Op}_{\varphi}(a)(x, hD_{x})$ is, modulo operators bounded in $L^2$ by ${\mathcal {O}}( h^\infty)$, independent of the choice of the function $\varphi$. To simplify notation 
 we drop writing $\varphi$.

Let us now consider a Riemannian manifold without boundary,  $M$. By partitions of unity we can define semi-classical pseudo-differential operators $a(x, hD_x)$ associated to symbols $a(x, \xi)\in \CIc (T^* M)$ 
\par
Now we consider a sequence $(u_n)$ bounded in $L^2(M)$. satisfying
\begin{equation}\label{eq.equation}
(-h_n^2\Delta -1)u_n = 0
\end{equation}  
Using~\eqref{eq.equation},
as in~\cite{GeLe93} (see also~\cite{BuIMRN}) we can prove the following
\begin{prop}
\label{prop2.1} There exist a subsequence $(n_{k})$ and a positive Radon  measure on $T^*M$, $\mu$ (a semi-classical measure for the sequence $(u_n)$), such that for any $a\in \CIc (T^* M)$
\begin{equation}
\label{eq2.8}\lim_{k\rightarrow + \infty}\left(\text{Op}(a)(x, h_{n_{k}}D_{x})f_{n_{k}}, f_{n_{k}}\right)_{ L^2(M}= \langle \mu, a(x, \xi)\rangle
\end{equation}
Furthermore this measure satisfies
\begin{enumerate}
\item The support of $\mu$ is included in the characteristic manifold:
\begin{equation}
\Sigma \stackrel{\rm{def}}{=} \{(x, \xi)\in T^*M; p(x,\xi)=\|\xi\|_x=1\}
\end{equation}
where  $\|\cdot \|_x$ is the norm for the metric at the point $x$
\item The measure $\mu$ is invariant by the bicharacteristic flow (the flow of the Hamilton vector field of $p$):
\begin{equation}
\label{eq.inv}H_p \mu =0
\end{equation}
\item For any $\varphi\in \CIc (M)$, 
\begin{equation}
\label{eq.lim} \lim_{n\rightarrow + \infty} \|\varphi u_n\|^2 = \langle \mu, |\varphi|^2\rangle
\end{equation}
\end{enumerate}
\end{prop}
The two first properties above are weak forms of the elliptic regularity and propagation of singularities results  whereas the last one states that there is no loss of $L^2$-mass at infinity in the $\xi$ variable.   
\section{Partially rectangular billiards}
\label{prb}
The following theorem is an easy conseauence of Proposition \ref{p:1}:
\begin{thm}
\label{t:1}
Let $ \Omega $ be a partially rectangular billiard with the rectangular
part $ R \subset \Omega $, $ \partial R = \Gamma_1 \cup \Gamma_2 $, a decomposition 
into parallel components satisfying $ \Gamma_2 \subset \partial \Omega $. 
Let $ \Delta $ be the Dirichlet or Neumann Laplacian on $\Omega $. Then
for {\em any} neighbourhood of $ \Gamma_1 $ in $ \Omega $, $ V $, there exists $ C $
such that
\begin{equation}
\label{eq:t1}
- \Delta u = \lambda u \ \Longrightarrow \  \int_V | u ( x ) |^2 dx 
\geq \frac1C \int_R | u ( x ) |^2 dx \,, 
\end{equation}
that is, no eigenfuction can concentrate in $ R $ and away from $ \Gamma_1$.
\end{thm}
\begin{proof}
Let us take $x,y$ as 
the coordinates on the stadium, so that $x$ parametrizes $ \Gamma_2 \subset \partial 
\Omega $ and $ y$, $ \Gamma_1 $, 
\[ R = [0,1]_{x}\times [0,a]_{y} \,.\]
Let $\chi\in \CIc((0,1))$ be equal to $1$ on $[\varepsilon, 1- \varepsilon]$. Then $\chi(x) u(x,y)$ is solution of
\begin{equation}
(\Delta-z)\chi u = [\Delta, \chi] u\text{ in $R$}
\end{equation} with the boundary conditions satisfied 
on $\partial R$. Applying Proposition~\ref{p:1}, we get 
\begin{equation}
\|\chi u\|_{L^2( R)}\leq C \left \| [ \Delta, \chi] u \|_{H^{-1}_{x}; L^2_{y}}+ 
\| u\rest _{\omega_{\varepsilon}}\|_{L^2( \omega_{\varepsilon})}\right)
\leq C' \| u\rest _{\omega_{\varepsilon}}\|_{L^2( \omega_{\varepsilon})} \,,
\end{equation}
where $\omega_{\varepsilon}$ is a neighbourhood of the support of $ \nabla\chi$.
Since a neighbourhood of $ \Gamma_1 $ in $ \Omega $ has to contain $ \omega_\varepsilon $
for some $ \varepsilon $, \eqref{eq:t1} follows.
\end{proof}

\section{Applications}
\label{app}

\begin{figure}[ht]
\includegraphics[width=10cm]{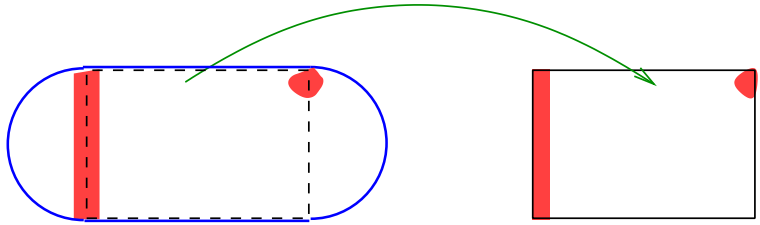}
\caption{Control regions in which eigenfunctions have positive mass
and the rectangular part for the Bunimovich stadium.}
\label{fig:bath}
\end{figure}

In \cite{BZ2} and \cite{BZ3} we used Proposition \ref{p:1} to prove that in the
case of the Bunimovich billiard shown in Fig.\ref{fig:bath}
 the states have nonvanishing 
density near the vertical boundaries of the rectangle. That follows from Theorem
\ref{t:1} which shows that we have to have positive density in the wings of the
billiard, and the propagation result (in the boundary case) based on the fact that
any diagonal controls a disc geometrically (see  \cite[Sect.6.1]{BZ2}; in fact 
we can use other control regions as shown in Fig.\ref{fig:bath}). 
Here we consider another case which accidentally generalizes a control theory 
result of Jaffard \cite{Ja}.

The Sinai billiard (see Fig.\ref{fig:sin}) 
is defined by removing a strictly convex open set, $ \mathcal O$, with 
a $ \CI $ boundary,
 from a flat torus, $ {\mathbb T}^2 \stackrel{\rm{def}}{=}
 \SP^1 \times \SP^1$:
\[ S  \stackrel{\rm{def}}{=} {\mathbb T}^2 \setminus {\mathcal O} \,.\]
Taking circles with different lengths 
might also possible but for simplicity we will restrict our attention
to a square torus.

\begin{thm}
\label{t:s}
Let $ V $ be any open neighbourhood of the convex boundary, $ \partial {\mathcal O} $, in 
a Sinai billiard, $ S $. If $ \Delta $ is the Dirichlet or Neumann
Laplace operator on $ S $ then there exists a constant, $ C = C ( V ) $, such that
\begin{equation}
\label{eq:ts}
- \Delta u = \lambda u \ \Longrightarrow \  \int_V | u ( x ) |^2 dx 
\geq \frac1C \int_S | u ( x ) |^2 dx \,. 
\end{equation}
\end{thm}
\begin{proof}
Suppose that the result is not true, that is, there exists a sequence of 
eigenfunctions $ u_n $, $ \| u _n \| = 1 $, with the corresponding 
eigenvalues $ \lambda_n \rightarrow \infty $,  such that $ \int_V | u_n ( x )|^2 dx
\rightarrow 0 $. We first observe that the only directions in the support of 
the corresponding semi-classical defect measure, $\mu$, have to be rational: the projection
of a trajectory with an irrational direction is dense on the torus and hence has
to encounter the obstacle $ \partial {\mathcal O} $ (and consequently $V$). The propagation result
recalled in Proposition \ref{prop2.1} gives a contradiction (remark that we apply this result as long as the trajectory does not encounter the obstacle and consequently we need only the {\em interior} propagation).

Hence let us assume that there exists a rational direction in the support of 
the measure which then contains the periodic trajectory in that direction.
As shown in Fig. \ref{f:3} we can find a maximal rectangular neighbourhood of 
the projection of that trajectory which avoids the obstacle: the sides parallel to the 
projection correspond to $ \Gamma_1 $ in Theorem \ref{t:1}.

\begin{figure}[ht]
$$\includegraphics[width=5.5cm]{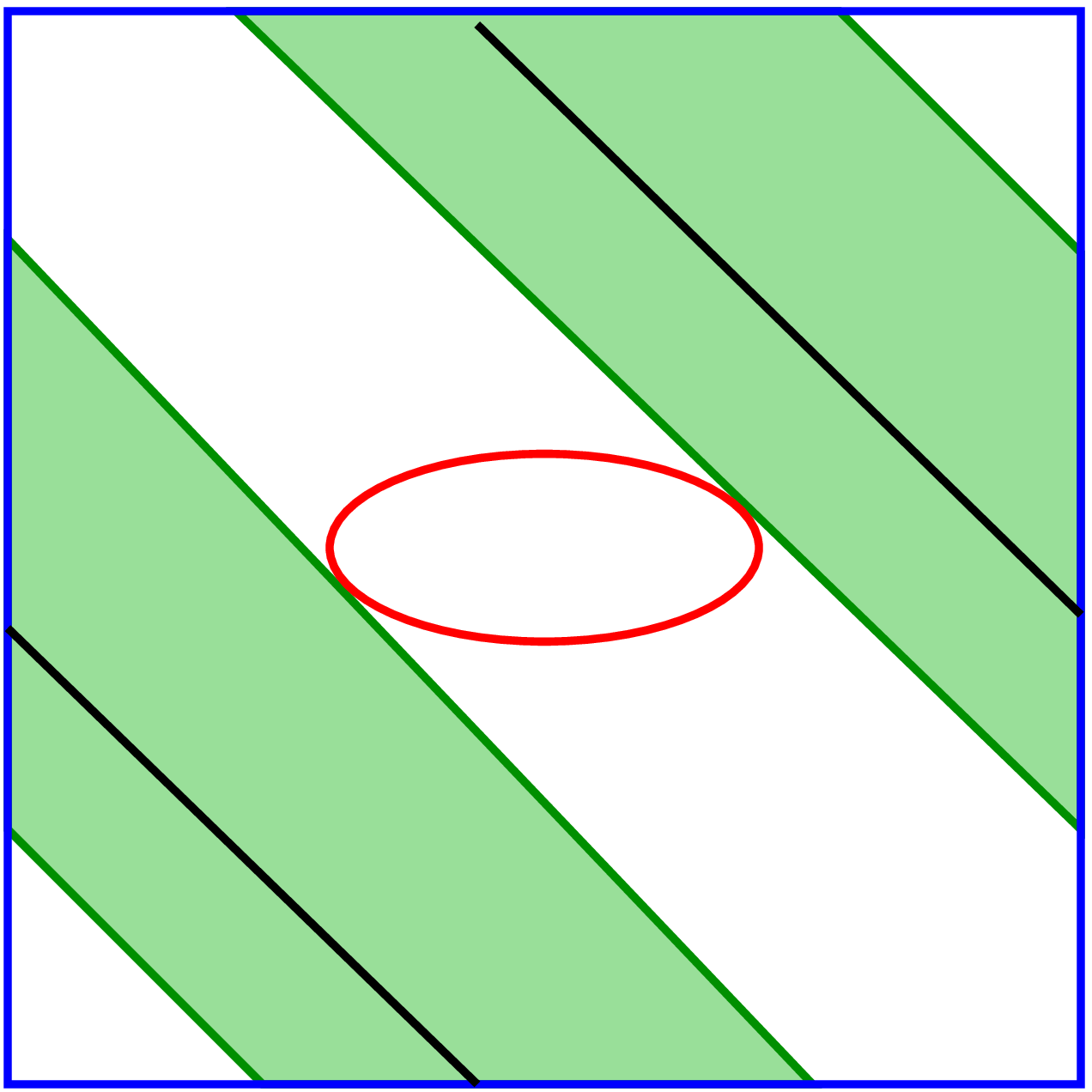}\qquad\includegraphics[height=5.5cm]{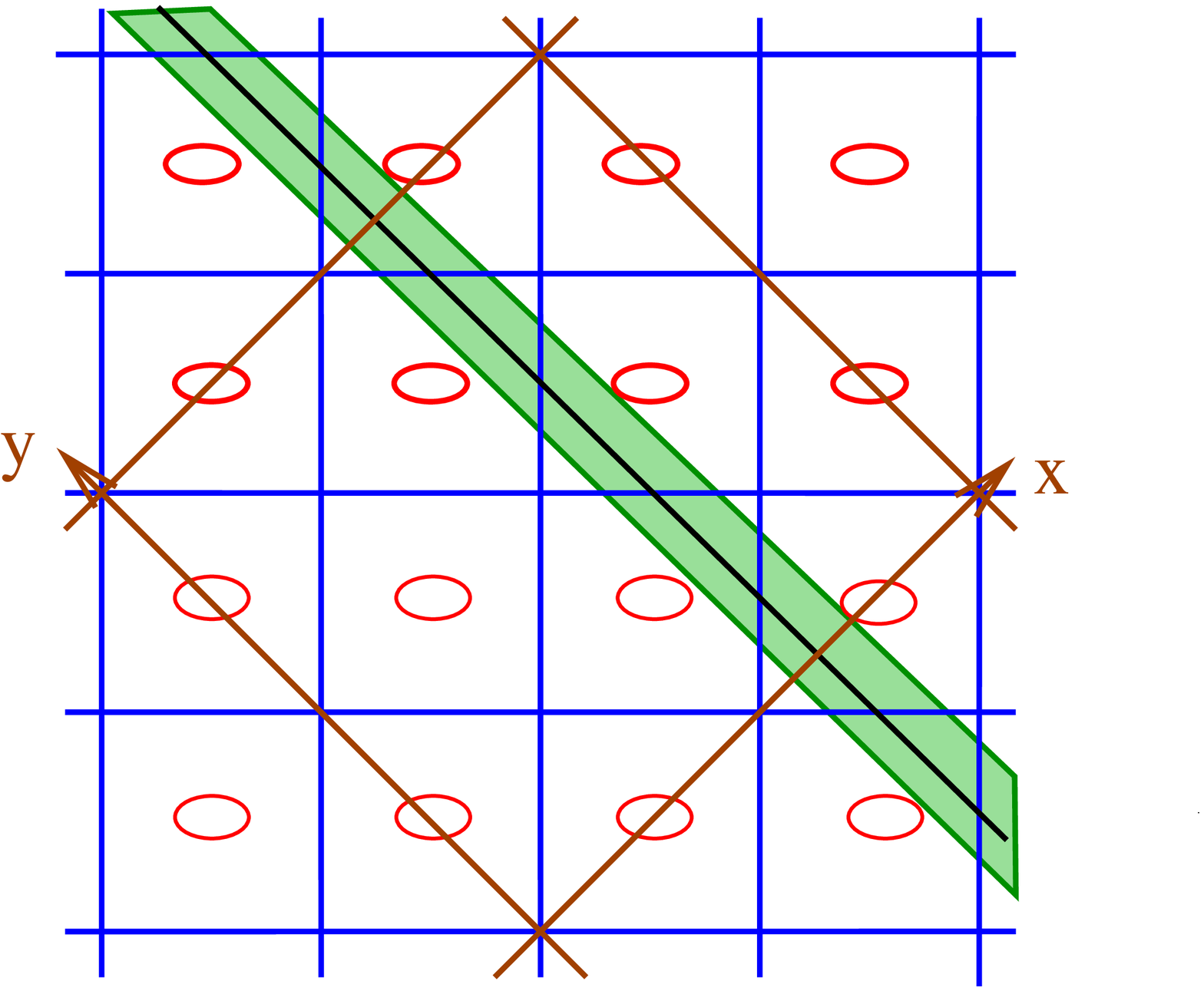}$$

\caption{A maximal rectangle in a rational direction, avoiding the obstacle. On the
right an explicit realization as a flat rectangle.}
\label{f:3}
\end{figure} 

The rectangle can be described as $ R = [0,a]_x \times [0, b]_y $ with the the 
$ y $ coordinate parametrizing the trajectory. 
Let $ u $ be an eigenfunction in our sequence and let $ \chi = \chi ( x) $ be 
a smooth function, supported in $ ( 0 , a ) $ and equal to one outside of a small
neighbourhood of the endpoints. Then $ \chi ( x ) u ( x, y ) $ is a function on 
$ R $ satisfying periodicity condition. Let $ E_\xi $ be a microlocal 
projection onto a neighbourhood of the $ R \times \{ \xi \} \subset T^* R $,
the semi-classical sense with $ h = 1/\sqrt \lambda $.
Let $ \Delta_R $ is the (periodic) Laplacian on $ R $. Using Fourier decomposition
we can arrange that $ [ \Delta_R , E_\xi ] = 0 $. Hence,
\[  ( - \Delta_R - \lambda ) E_\xi \chi u = [ \Delta_R , E_\xi \chi ] \tilde \chi u 
=E_\xi  [ \Delta_R , \chi ] \widetilde E_\xi \tilde \chi u + {\mathcal O} ( \lambda
^{- \infty } ) 
\,, \ \  \| u \| = 1 \,, 
\]
where $ \tilde \chi $ has the same properties as $ \chi $ and is equal to 
one on the support of $ \chi $, and similarly for $ \widetilde E_\xi $.
As in the proof of Theorem \ref{t:1} and using that $E_\xi$ is continuous on $H^{-1}_x; L^2_y$, we now see
that 
\begin{equation}
\label{eq:mic}
 \| E_\xi \chi u \| \leq C \int_\omega | \widetilde E_\xi \tilde \chi u|^2 
+ {\mathcal O} ( \lambda^{-\infty } ) \,,
\end{equation}
where $ \omega $ is a neighbourhood of $ \nabla \chi $ (we are using here
the calculus of semi-classical pseudo-differential operators). 
Since the semi-classical defect measure of $ E_\xi \chi u $ (which is $|E_\xi \chi |^2\times \mu$) was assumed to be
non-zero \eqref{eq:mic} shows that the measure of $ \widetilde E_\chi \tilde \chi u 
\rest_\omega $
is non zero and consequently there is a point in the intersection of the supports of $\mu$ and $\widetilde e_\xi \tilde \chi$. But $\mu$ is invariant
by the flow (as long as it does not intersect the obstacle) and hence, once we choose all the cut-offs above very close to the
boundary of $ R$, its support can be made intersect any neighbourhood of $ \partial {\mathcal O} $.
\end{proof}
\begin{rem}
In the proof above the smoothness, the convexity, and even the connectivity
  of the obstacle played no role (and we could take $\Theta= \emptyset$ provided that $V\neq \emptyset$). Consequently, the result holds for any obstacle (sufficiently smooth in the 
case of Neumann boundary conditions) 
and consequently to the special case of pseudointegrable billiards (see for instance
\cite{Bog} for motivation and description). 
By an elementary reflection principle, the result also holds for an obstacle inside a square with Dirichlet or Neumann conditions on the boundary of the square.
\end{rem}
\begin{rem}
The proof above gives in fact the following estimate for any open neighbourhood of the obstacle:
\begin{equation}
\begin{gathered}
\exists C; \forall u,\, f \in L^2(S) \text{ solutions of } (-\Delta + \lambda)u =f, \qquad u \rest_{\partial S} =0\\
\|u\|_{L^2(S)}\leq C \left( \|f\|_{L^2(S)}+ \|u \11 _V\|_{L^2(V)}\right)
\end{gathered} 
\end{equation}
and according to~\cite[Theorem 4]{BZ2}, this implies that the Schr{\"o}dinger equation in $S$ is exactly controlable by $V$ in finite time. In fact, by working on the time evolution equation, we could strenghten this result allowing an arbitrarily small time. This latter result was previously known~\cite{Ja} for the particular case $\Theta= \emptyset$ ($S= \mathbb{T}^2$) but the proof was 
based on subtle results about Fourier series \cite{Ka}.
\end{rem}
\def\cprime{$'$}
\begin{rem}
As shown in \cite[Theorem 2$'$]{BZ2}, the results of Ikawa and G\'erard on 
scattering by two convex obstacles (see \cite{BZ2} and references given there)
give an estimate on the maximal concentration of an eigenfunction (or a 
quasimode) on a closed orbit in a Sinai billiard. Let $ \chi \in \CI(S;[0,1])  $  be 
supported in  a small neighbourhood of a closed transversally reflecting 
orbit. Then  for any family $ ( - \Delta - \lambda ) u_\lambda = {\mathcal O}
( \lambda^{-\infty } ) $, $ \| u_\lambda \| = 1 $, 
\[ C \int_S  | u ( x ) |^2 ( 1 - \chi ( x ) ) dx \geq \frac{1}{ \log \lambda } \,, \]
that is a concentration on a closed trajectory, if at all possible, has to 
be very weak.
\end{rem}

\end{document}